\documentclass[10pt,twoside]{article}
\usepackage{amsmath,amssymb,amsthm}
\usepackage{hyperref}

\usepackage[hyperref,backref,natbib,ibidtracker=false]{biblatex}
\bibliography{secchi}

\usepackage{latexsym}
\newtheorem{theorem}{Theorem}[section]
\newtheorem{definition}{Definition}[section]
\newtheorem{lemma}{Lemma}[section]
\newtheorem{proposition}{Proposition}[section]
\newtheorem{remark}{Remark}[section]

\def\proof{\mbox {\it Proof.~}}

\let\ge=\varepsilon

\makeatletter\def\theequation{\arabic{section}.\arabic{equation}}\makeatother
\begin{document}
\title{
{\bf\Large  The Brezis--Nirenberg problem for the H\'{e}non equation: ground state solutions}
\author{{\bf\large Simone Secchi}\footnote{Partially supported by PRIN 2009 ``Teoria dei punti critici
e metodi perturbativi per equazioni differenziali nonlineari''.}\\
\hspace{2mm}
{\it \small Dipartimento di Matematica ed Applicazioni} \\
{\it \small Universit\`a di Milano--Bicocca}\\
{\it\small Via R. Cozzi 53, I-20125 Milano (Italy)}\\
{\it \small e-mail: Simone.Secchi@unimib.it}}
}

\maketitle
\begin{center}
{\bf\small Abstract}

\vspace{3mm}
\hspace{.05in}\parbox{4.5in}
{{\small This work is devoted to the Dirichlet problem for the 
       equation \(-\Delta u = \lambda u + |x|^\alpha |u|^{2^*-2} u\) in the unit ball of $\mathbb{R}^N$. 
       We assume that $\lambda$ is bigger than the first eigenvalues of the laplacian, and we prove that
there exists a solution provided $\alpha$ is small enough. This solution has a variational characterization as a ground state.}}
\end{center}


\section{\bf Introduction}
\def\theequation{1.\arabic{equation}}\makeatother
\setcounter{equation}{0}
This short note is devoted to the Dirichlet problem
    \begin{equation} \label{eq:1}
    \begin{cases}
    -\Delta u = \lambda u + |x|^\alpha |u|^{2^*-2}u &\hbox{in $\Omega$}\\
    u = 0 &\hbox{on $\partial \Omega$}
    \end{cases}
    \end{equation}
where $\Omega$ is the unit ball of $\mathbb{R}^N$, $\lambda$ is bigger than $\lambda_1=\lambda_1(-\Delta)$, the first Dirichlet eigenvalue of $-\Delta$, and $\alpha$ is a positive parameter. The exponent $2^*$ is a shorthand for the Sobolev critical exponent $2N/(N-2)$. We will assume throughout that $N \geq 3$.

This problem is a generalization of the celebrated Brezis--Nirenberg problem, see \cite{MR709644} and \cite{MR829403,MR831041,MR1441856,MR2122698} for more general and/or recent existence results. When $\alpha \neq 0$, our equation is reminiscent of the H\'{e}non equation
\[
-\Delta u = |x|^\alpha |u|^{p-2}u,
\]
which has been studied deeply in recent times. Most papers deal with the subcritical case $p<(N+2)/(N-2)$, and focus on the behavior of solutions as $\alpha \to +\infty$ or $p \to (N+2)/(N-2)$. We refer to \cite{MR1918755,MR2507301,MR2541412,MR1963460,MR2507301} for more information. As far as we know, the Brezis--Nirenberg problem for the critical H\'{e}non equation has been studied only in  \cite{LY}, where the authors prove that there always exists a solution to problem (\ref{eq:1}), provided $N \geq 7$ and $\alpha$ is small enough.

In the next sections we will show that solutions exist whenever $N \geq 5$ and $\alpha$ is small; in addition, we will find them as ground-state solutions, in a sense that will be made precise in a moment. We can therefore remove the (technical) restriction on the space dimension, and also provide more information about solutions.
We will borrow many ideas from the recent papers~\cite{Pankov} and~\cite{MR2553063}, although the presence of the increasing weight $|\cdot|^\alpha$ has to be dealt with carefully.
Our main result is the following theorem.
\begin{theorem} \label{th:main}
Denote by $\lambda_1 < \lambda_2 \leq \lambda_3 \ldots \leq \lambda_m\leq \ldots$ the Dirichlet eigenvalues of the laplacian, and assume that $N \geq 5$. If $\lambda_m \leq \lambda < \lambda_{m+1}$ for some \(m \in \mathbb{N}\), then, for every $\alpha>0$ sufficiently small,
there exists (at least) a ground-state solution to problem (\ref{eq:1}).
\end{theorem}
For the precise definition of ground-state solutions, we refer to Definition \ref{def:gs} below.

As a consequence of well-known results in bifurcation theory for potential operators (we refer to Theorem 6.1 of \cite{Amb}), it is rather easy to prove that each eigenvalue~$\lambda_m$ is a bifurcation point for problem (\ref{eq:1}): this is the reason why many papers focused on the case $\lambda \notin \sigma (-\Delta)$. We propose a variational approach that also covers the case $\lambda = \lambda_m \in \sigma(-\Delta)$.

\section{\bf A variational framework for ground-state solutions}
\def\theequation{2.\arabic{equation}}\makeatother
\setcounter{equation}{0}
We will work in the Hilbert space $H=H_0^1(\Omega)$ endowed with the Dirichlet inner product
\[
\langle u , v \rangle = \int_\Omega \nabla u \cdot \nabla v 
\]
and the induced norm $\| \cdot \|$. We will assume that, for some~$m \in \mathbb{N}$,
\[
\lambda_m \leq \lambda < \lambda_{m+1},
\]
as stated in Theorem \ref{th:main}. We denote by $\{e_j\}_j$ the eigenfunctions associated to \(\{\lambda_j\}_j\). By assumption, we are led to the decomposition
\[
H = Z \oplus Y,
\]
where $Z$ is the subspace of $H$ spanned by the first $m$ eigenfunctions $e_1$,\ldots,$e_m$ and $Y = Z^\perp$. There is a standard identification of solutions to (\ref{eq:1}) with the critical points of the functional $\varphi \colon H \to \mathbb{R}$ defined by the formula
\begin{equation}
\varphi(u) = \frac{1}{2} \int_\Omega \left( |\nabla u|^2 - \lambda |u|^2 \right) - \frac{1}{2^*} \int_\Omega |x|^\alpha |u|^{2^*}.
\end{equation}
In order to find \emph{ground state} solutions of (\ref{eq:1}), we introduce (see \cite{Pankov}) a sub-manifold of $H$,
\begin{equation}
 \mathcal{N} = \left\lbrace u \in H \setminus \{0\} \mid \langle \nabla \varphi (u),u \rangle =0, \ \nabla \varphi (u) \in Y \right\rbrace
\end{equation}
\begin{remark}
The set \(\mathcal{N}\) is the intersection of the standard \emph{Nehari manifold} 
\[
\left\{u \in H \mid \langle \nabla \varphi(u),u\rangle =0 \right\}
\]
with the pre-image \(\left(\nabla \varphi \right)^{-1}(Y)\). Much more general cases of Nehari-like manifolds and natural constraints are studied in \cite{NV}.
\end{remark}

\begin{proposition} \label{prop:2.1}
 The set $\mathcal{N}$ is a $C^1$ submanifold of $H$, of codimension $m+1$. Moreover, $\mathcal{N}$ is a natural
contraint for $\varphi$: every critical point of the restriction $\varphi_{|\mathcal{N}}$ is a free critical point of $\varphi$.
\end{proposition}
\noindent\proof
 We borrow the proof from \cite{MR2553063}. Consider the map $F \colon H \setminus \{0\} \to \mathbb{R} \times Z$, defined by the formula
\[
 F(u)=(\langle \nabla \varphi(u),u\rangle,Q\nabla \varphi(u)),
\]
where $Q$ is the orthogonal projection of $H$ onto $Z$; then $\mathcal{N} = F^{-1}(0)$. On the cartesian product $\mathbb{R} \times Z \simeq \mathbb{R}^{m+1}$ we put the inner product
\[
 (t_1,z_1) \cdot (t_2,z_2) = t_1 t_2 + \langle z_1,z_2 \rangle.
\]
We claim that
\[
 \left( DF(u)(tu+z) \right) \cdot (t,z) <0
\]
for any $(t,z) \in \mathbb{R} \times Z$, $(t,z) \neq (0,0)$. It is elementary to realize that this claim completes the proof of the first part of our Proposition. Fix $(t,z) \neq (0,0)$, and remark that
\[
 \langle \nabla \varphi(u),u \rangle = \langle \nabla \varphi (u),z \rangle =0
\]
implies
\begin{multline*}
 \left( DF(u)(tu+z) \right) \cdot (t,z) \\
= t D^2 \varphi(u)(tu+z,u)+t \langle \nabla \varphi(u),tu+z \rangle + D^2 \varphi(u)(tu+z,z)\\
= D^2 \varphi(u)(tu+z,tu+z) - t \langle \nabla \varphi(u),tu+2z \rangle \\
= \int_\Omega |\nabla z|^2 - \lambda |z|^2 \, dx - \int_\Omega \left( (2^*-1)(tu+z)^2 - tu(tu+2z)\right) |u|^{2^*-2} |x|^\alpha  dx \\
= \int_\Omega |\nabla z|^2 - \lambda |z|^2 \, dx \\
- \int_\Omega \left( (2^*-2)t^2 u^2 + 2 (2^*-2)tzu + (2^*-1)z^2 \right) |u|^{2^*-2} |x|^\alpha \, dx.
\end{multline*}
As a quadratic form in $(t,z)$, the integral 
\[
\int_\Omega \left( (2^*-2)t^2 u^2 + 2 (2^*-2)tzu + (2^*-1)z^2 \right) |u|^{2^*-2} |x|^\alpha \, dx
\]
 is positive definite whenever $u(x) \neq 0$. By the assumption $\lambda_m \leq \lambda < \lambda_{m+1}$, the quadratic form $\int_\Omega |\nabla z|^2 - \lambda |z|^2$ is negative semidefinite.
If $\int_\Omega |\nabla z|^2 - \lambda |z|^2 < 0$, the claim is proved. If $\int_\Omega |\nabla z|^2 - \lambda |z|^2=0$,
either $\lambda = \lambda_m$ and $z$ is an eigenfunction, or $z=0$. By assumption, $t \neq 0$ if $z=0$; moreover, $z \neq 0$ implies
$z \neq 0$ almost everywhere. In both cases, the claim follows easily.

Finally, we need to check that $u \in H$ is a critical point of $\varphi$ if and only if $u \in \mathcal{N}$ and $D\varphi(u)$ vanishes on
the tangent space $T_u \mathcal{N}$. The necessary condition is trivial; on the contrary, assuming that $D\varphi(u) =0$ on $T_u \mathcal{N}$ and $u \in \mathcal{N}$,
we deduce that $D\varphi(u)$ also vanishes on $\mathbb{R}u \oplus Z$. But we have just proved that $\mathbb{R}u \oplus Z$ is transversal
to $T_u \mathcal{N}$, and we conclude.

\begin{remark}
 The previous Proposition states that $DF(u)$ is a surjective map at every $u \in F^{-1}(0)\setminus \{0\}$. But the additional
information that $\left( DF(u)(tu+z) \right) \cdot (t,z)$ is \emph{negative} will be useful later on.
\end{remark}

Since $\mathcal{N}$ contains every critical point of $\varphi$, the following terminology is rather natural.
\begin{definition} \label{def:gs}
 A \emph{ground state solution} $u$ to (\ref{eq:1}) is any element of $\mathcal{N}$ such that $D\varphi(u)$ vanishes on $T_u \mathcal{N}$ and
$\varphi (u)=c$, where the level $c$ is defined by
\begin{equation} \label{eq:c}
 c = \inf_{\mathcal{N}} \varphi 
\end{equation}
\end{definition}

The arguments of \cite{MR2557725}, which hold true under general assumptions, guarantee that for every $v \in Y \setminus\{0\}$
there exists a unique couple $(f(v),g(v)) \in (0,+\infty) \times Z$ such that $F(f(v)v+g(v))=0$. Moreover $f(\cdot)$ and $g(\cdot)$ 
are continuous maps, and
\[
 \varphi(f(v)v+g(v))= \max_{\substack{t > 0 \\ w \in Z}} \varphi (tv+w).
\]

It follows easily from the definition of $f$ and $g$ that
\[
 c = \inf_{\substack{v \in Y \\ v \neq 0}} \varphi(f(v)v+g(v)) = \inf_{\substack{v \in Y \\ v \neq 0}} \max_{\substack{t > 0 \\ w \in Z}}
\varphi (tv+w).
\]

\section{\bf Existence of ground state solutions}
\def\theequation{3.\arabic{equation}}\makeatother
\setcounter{equation}{0}
The existence of a ground state solution to (\ref{eq:1}) will be proved by a compactness argument. Since (\ref{eq:1})
contains the critical exponent, it is natural to expect compactness of minimizing sequences (for $c$) below some 
energy level related to Sobolev's best constant $S$. Recall that
\[
 S = \inf_{\substack{u \in H \\ u \neq 0}} \frac{\int_\Omega |\nabla u|^2 \, dx}{\left( \int_\Omega |u|^{2^*}\, dx \right)^{2/2^*}},
\]
and this numer is actually indipendent of the domain~$\Omega$. A simple exercise in sophomore calculus proves the next lemma, stated in \cite{MR2553063}.
\begin{lemma} \label{lemma:calculus}
 If $A>0$ and $B>0$, then
\[
 \max_{t>0} \left( \frac{1}{2}At^2 - \frac{1}{2^*}Bt^{2^*} \right) = \frac{1}{N} \left( \frac{A}{B^{2/2^*}} \right)^{N/2}.
\]
\end{lemma}
We now come to the main compactness result about the variational problem (\ref{eq:c}).
\begin{proposition} \label{prop:3.2}
 Suppose that
\begin{equation} \label{eq:5}
 c < \frac{1}{N} S^{N/2}.
\end{equation}
Then there exists $v \in Y \setminus \{0\}$ such that
\[
 \max_{\substack{t >0 \\ w \in Z}} \varphi(tv+w) = \varphi(f(v)v+g(v))=c.
\]
\end{proposition}
\noindent\proof
Take any sequence $\{v_n\}_n$ in $Y \setminus \{0\}$ such that $\|v_n\|=1$ and
\begin{equation} \label{eq:6}
\max_{\substack{t>0 \\ w \in Z}} \varphi(tv_n+w) \to c.
\end{equation}
Without loss of generality, we can assume that $v_n \to v$ weakly in $H$, strongly in $L^2(\Omega)$ and point-wise almost everywhere.
Writing
\begin{align*}
A &= \lim_{n \to +\infty} \int_\Omega |\nabla (v_n-v)|^2 \, dx \\
B_\alpha &= \lim_{n \to +\infty} \int_\Omega |v_n-v|^{2^*} |x|^\alpha \, dx
\end{align*}
and using the Brezis--Nirenberg lemma, we exploit (\ref{eq:6}) to get
\begin{equation} \label{eq:7}
\varphi (tv+w) + \frac{1}{2} A t^2 - \frac{1}{2^*} B_\alpha t^{2^*} \leq c.
\end{equation}
We now distinguish several possibilities. If $v=0$ and $B_\alpha =0$, from the assumption $\|v_n\|=1$ we deduce $A=1$. Hence $t^2 \leq 2c$ for every $t>0$, a contradiction.

Assume now $B_\alpha \neq 0$. From the Sobolev inequality and the trivial remark that $|x|^\alpha <1$ in $\Omega$, we get
\begin{equation} \label{eq:sob}
\frac{1}{N}S^{N/2} \leq \frac{1}{N} \left( \frac{A}{B_0^{2/2^*}} \right)^{\frac{N}{2}} \leq \frac{1}{N} \left( \frac{A}{B_\alpha^{2/2^*}} \right)^{\frac{N}{2}} = \max_{t>0} \left( \frac{1}{2}A t^2 - \frac{1}{2^*}B_\alpha t^{2^*} \right).
\end{equation}
If $v=0$, we conclude that
\[
\frac{1}{N}S^{N/2} \leq c < \frac{1}{N}S^{N/2},
\]
and thus $v \neq 0$. Call $h = g(v)/f(v)$. It follows from the definition of the level $c$ that
\begin{multline} \label{eq:8}
c \leq \varphi(f(v)(v+h)) = \max_{t>0} \varphi (t(v+h)) \\
= \frac{1}{N} \frac{\int_\Omega |\nabla v|^2 + |\nabla h|^2 - \lambda (v^2+h^2) \, dx}{\left( \int_\Omega |v+h|^{2^*} |x|^\alpha \, dx \right)^{2/2^*}}.
\end{multline}
From (\ref{eq:7}),
\begin{multline} \label{eq:9}
\max_{t>0} \left( \varphi(t(v+h)) + \frac{1}{2} A t^2 - \frac{1}{2^*} B_\alpha t^{2^*} \right) \\
= \frac{1}{N} \frac{A+\int_\Omega |\nabla v|^2 + |\nabla h|^2 - \lambda (v^2+h^2) \, dx}{\left(B_\alpha+ \int_\Omega |v+h|^{2^*} |x|^\alpha \, dx \right)^{2/2^*}} \leq c.
\end{multline}
Putting together (\ref{eq:5}), (\ref{eq:sob}), (\ref{eq:8}) and (\ref{eq:9}) we can write
\begin{multline}
(Nc)^{2/N} \left( B_\alpha + \int_\Omega |v+h|^{2^*}|x|^\alpha \, dx \right)^{2/2^*} \\
< (Nc)^{2/N} \left( B_\alpha^{2/2^*} + \left( \int_\Omega |v+h|^{2^*} |x|^\alpha \, dx \right)^{2/2^*} \right) \\
< A + \int_\Omega |\nabla v|^2 + |\nabla h|^2 - \lambda (v^2+h^2) \, dx \\
\leq (Nc)^{2/N} \left( B_\alpha + \int_\Omega |v+h|^{2^*} |x|^\alpha \, dx \right)^{2/2^*},
\end{multline}
a contradiction. Therefore $B\sb\alpha=0$ and (\ref{eq:7}) yields
\[
c \leq \varphi(f(v)v+g(v)) \leq c.
\]
\begin{remark}
It was proved in \cite{LY}, mimicking the ideas contained in \cite{MR709644}, that $\varphi$ satisfies the Palais--Smale condition below the threshold $S^{N/2}/N$. The same result could also be proved by slightly adapting the arguments of \cite{MR2142066}.
\end{remark}

The subspace $Z$ has a kind of unique continuation property, as proved in \cite[Lemma 3.3]{MR2553063}.

\begin{lemma}\label{lemma:3.4}
If $w \in Z$ vanishes on some open subset $\omega \neq \emptyset$ of $\Omega$, then $w=0$ everywhere.
\end{lemma}

The next step is to check that the level $c$ defined in (\ref{eq:c}) satisfies inequality (\ref{eq:5}). We proceed in several steps.

Consider $\ell \ll 1$, a parameter that will tend to zero at a slower rate than $\ge$: $\ge/\ell \to 0$. As $\ell \to 0$, the point
\[
x_\ell = (1-\ell,0,\ldots,0)
\]
approaches the boundary of $\Omega$. We pick a test function $\xi=\xi_\ell\in C_0^\infty(\Omega)$ spiked at $x_\ell$:
\begin{equation*}
\xi(x) = 
\begin{cases}
1, &\hbox{if $x \in B(x_\ell,\ell/2)$}\\
0, &\hbox{if $x \notin B(x_\ell,\ell)$},
\end{cases}
\end{equation*}
and such that $|\nabla \xi_\ell| \leq C / \ell$.
It is well known that
the \emph{instanton} $U_{\ge}$ defined by the formula
\begin{equation*}
U_{\ge,\ell} (x) = \left(N(N-2) \right)^{\frac{N-2}{4}} \frac{\ge^{\frac{N-2}{2}}}{\left(\ge^2 +|x-x_\ell|^2\right)^{\frac{N-2}{2}}}
\end{equation*}
is the optimal function for the Sobolev inequality in 
\[
D^{1,2}(\mathbb{R}^N) = \left\{u \in L^{2^*}(\mathbb{R}^N) \mid \nabla u \in L^2(\mathbb{R}^N) \right\}.
\]
Call now
\[
u_{\ge,\ell}(x) = \xi_\ell(x)U_{\ge,\ell}(x).
\]
Up to a constant that we can neglect in the following estimates, we can pretend that
\[
u_{\ge,\ell}(x) = \xi_\ell(x) \frac{\ge^{\frac{N-2}{2}}}{\left(\ge^2 +|x-x_\ell|^2\right)^{\frac{N-2}{2}}}.
\]
Reasoning as in \cite{MR709644}, we can estimates
\begin{equation} \label{est:1}
\int_\Omega |\nabla u_{\ge,\ell}|^2 \, dx = S^{\frac{N}{2}} + \hbox{h.o.t.},
\end{equation}
where h.o.t denotes higher order terms like
\begin{multline*}
\ge^{N-2} \int_{B(x_\ell,\ell) \setminus B(x_\ell,\ell/2)} \frac{|\nabla \xi_\ell (x)|^2}{\left(\ge^2+|x-x_\ell|^2 \right)^{N-2}} dx \\
\leq C \frac{\ge^{N-2}}{\ell^2} \int_{B(x_\ell,\ell) \setminus B(x_\ell,\ell/2)}  \frac{dx}{\left(\ge^2+|x-x_\ell|^2 \right)^{N-2}} \\
= C \frac{\ge^{N-2}}{\ell^2}  \int_{B(0,\ell) \setminus B(0,\ell/2)}  \frac{dy}{\left(\ge^2+|y|^2 \right)^{N-2}} \\
\leq C \frac{\ge^{N-2}}{\ell^2} \int_{B(0,\ell) \setminus B(0,\ell/2)} \frac{dy}{\left(\ge^2+\frac{\ell^2}{16}\right)^{N-2}} \\
= C \frac{\ge^{N-2}}{\ell^2} \frac{\ell^N}{\ell^{2(N-2)}} = C \left( \frac{\ge}{\ell} \right)^{N-2}.
\end{multline*}
As a consequence,
\begin{equation} \label{est:2}
\int_\Omega |\nabla u_{\ge,\ell}|^2 \, dx = S^{\frac{N}{2}} + O\left( \left( \frac{\ge}{\ell} \right)^{N-2} \right).
\end{equation}
Similarly,
\begin{multline*}
\int_\Omega |u_{\ge,\ell}|^{2^*} = \ge^N \int_\Omega \frac{|\xi_{\ell}(x)|^{2^*}}{\left(\ge^2+|x-x_\ell|^2 \right)^N} \\
= \ge^N \int_\Omega \frac{|\xi_\ell (x)|^{2^*} -1}{\left(\ge^2+|x-x_\ell|^2 \right)^N} \, dx + \ge^N \int_\Omega \frac{dx}{\left(\ge^2+|x-x_\ell|^2 \right)^N} \\
= \int_{\mathbb{R}^N} |U_{\ge,\ell}|^{2^*} + O(\ge^N) + \ge^N \int_\Omega \frac{|\xi_\ell (x)|^{2^*} -1}{\left(\ge^2+|x-x_\ell|^2 \right)^N} \, dx.
\end{multline*}
But
\begin{multline*}
\ge^N \int_\Omega \frac{|\xi_\ell (x)|^{2^*} -1}{\left(\ge^2+|x-x_\ell|^2 \right)^N} \, dx = 
\ge^N \int_{\Omega \setminus B(x_\ell,\ell/2)} \frac{|\xi_\ell (x)|^{2^*} -1}{\left(\ge^2+|x-x_\ell|^2 \right)^N} \, dx \\
\leq C \frac{\ge^N}{\ell^{2N}} \ell^N = C \left( \frac{\ge}{\ell} \right)^N.
\end{multline*}
We conclude that
\begin{equation} \label{est:3}
\int_\Omega |u_{\ge,\ell}|^{2^*} = \int_{\mathbb{R}^N} |U_{\ge,\ell}|^{2^*} + O(\ge^N) + O\left( \left( \frac{\ge}{\ell} \right)^N \right).
\end{equation}
The $L^2$-norm is slightly more involved:
\begin{multline*}
\int_\Omega |u_{\ge,\ell}|^2 = \ge^{N-2} \int_{B(x_\ell,\ell)} \frac{dx}{\left(\ge^2+|x-x_\ell|^2 \right)^{N-2}} \\
+ \ge^{N-2} \int_{B(x_\ell,\ell) \setminus B(x_\ell,\ell/2)} \frac{\xi_\ell(x)^2}{\left(\ge^2+|x-x_\ell|^2 \right)^{N-2}}.
\end{multline*}
Now,
\begin{multline*}
\ge^{N-2} \int_{B(x_\ell,\ell) \setminus B(x_\ell,\ell/2)} \frac{\xi_\ell(x)^2}{\left(\ge^2+|x-x_\ell|^2 \right)^{N-2}}
\leq C \ge^{N-2} \int_{\ell/2}^\ell r^{3-N}\, dr \\
= C \frac{\ge^N-2}{\ell^{N-4}}.
\end{multline*}
On the other hand,
\begin{multline*}
\ge^{N-2} \int_{B(x_\ell,\ell)} \frac{dx}{\left(\ge^2+|x-x_\ell|^2 \right)^{N-2}} \geq C \ge^{N-2} \int_{B(x_\ell,\ge)}
\frac{dx}{(2\ge^2)^{N-2}} + \\
C \ge^{N-2} \int_{B(x_\ell,\ell)\setminus B(x_\ell,\ge)} \frac{dx}{(2|x-x_\ell|^2)^{N-2}} = C\ge^2 + C \ge^{N-2}
\int_\ge^\ell r^{3-N}\, dr \\
= C\ge^2 + O\left( \frac{\ge^{N-2}}{\ell^{N-4}}\right).
\end{multline*}
We are now ready to estimate
\begin{multline} \label{eq:15}
\frac{\int_\Omega |\nabla u_{\ge,\ell}|^2 - \lambda |u_{\ge,\ell}|^2}{\left(\int_\Omega |x|^\alpha |u_{\ge,\ell}|^{2^*} \right)^{2/2^*}}
\leq \frac{S^{\frac{N}{2}} + O\left(\left(\frac{\ge}{\ell}\right)^{N-2}\right)-\lambda \left( C\ge^2 + O\left( \frac{\ge^{N-2}}{\ell^{N-4}} \right) \right)}{(1-2\ell)^{\alpha \frac{2}{2^*}} \left( S^{\frac{N}{2}} + O ((\ge/\ell)^{N})\right)^{2/2^*}} \\
= \frac{1}{(1-2\ell)^{\frac{2\alpha}{2^*}}} \left( S - C \ge^2 + O\left( \left(\frac{\ge}{\ell} \right)^{N-2}  \right)   \right)
\end{multline}

\begin{proposition} \label{prop:3.5}
There results 
\[
c< \frac{1}{N} S^{\frac{N}{2}}.
\]
\end{proposition}
\noindent\proof
We will check that
\begin{equation} \label{eq:16}
\max_{\substack{t>0 \\ w \in Z}} \varphi(tu_{\ge,\ell}+w) < \frac{1}{N}S^{N/2}.
\end{equation}
Setting $\omega = \Omega \setminus \operatorname{supp} \xi_{\ell}$, Lemma \ref{lemma:3.4} implies that $w \mapsto \|w\|_{L^{2^*}(\omega)}$ defines a norm on the subspace $Z$. Since $\dim Z = m < +\infty$, all norms on $Z$ are equivalent: we will use this remark tacitly in the sequel. 

We choose $\ell=\sqrt[4]{\ge}$, and write $u_{\ge}$ instead of $u_{\ge,\ell}$.
By elementary convexity and recalling that $0 \notin \omega$ so that $|x|^\alpha$ is bounded away from zero as $x \in \omega$, for every $t>0$ and every $w \in Z$, we can estimate
\begin{multline*}
\int_\Omega |tu_{\ge}+w|^{2^*}|x|^\alpha \, dx = \int_{\Omega\setminus\omega} |tu_{\ge}+w|^{2^*}|x|^\alpha \, dx + \int_\omega |w|^{2^*} |x|^\alpha \, dx \\
\geq t^{2^*} \int_\Omega |u_{\ge}|^{2^*} |x|^\alpha \, dx + 2^* t^{2^*-1} \int_\Omega |u_{\ge}|^{2^*-1} w |x|^\alpha \, dx + 2^* C_1 \|w\|^{2^*}.
\end{multline*}
It follows that
\begin{multline} \label{eq:17}
\varphi(t u_{\ge}+w) \leq \varphi(t u_{\ge}) + t \int_{\Omega} \nabla u_{\ge} \cdot \nabla w - \lambda u_{\ge,\ell}w \, dx \\
+ \frac{1}{2} \int_{\Omega} |\nabla w|^2 - \lambda |w|^2 \, dx - t^{2^*-1} \int_{\Omega} |u_{\ge}|^{2^*-1} w |x|^\alpha \, dx - C_1 \|w\|^{2^*}.
\end{multline}
In particular, we can write
\[
\varphi (t u_{\ge}+w) \leq A \left( t^2 + t \|w\| + t^{2^*-1}\|w\| \right) - B \left( t^{2^*}+\|w\|^{2^*} \right)
\]
for suitable constants $A>0$ and $B>0$. Hence there exists a number $R>0$ such that, for $\ge$ and $\ell$ small, $t>R$ and $w \in Z$ there holds $\varphi(tu_{\ge}+w) \leq 0$. On the other hand, whenever $t \leq R$,
\[
\varphi(tu_{\ge}+w) \leq \varphi(t u_{\ge}) + O(\ge^{\frac{N-2}{2}}) \|w\| - C_1 \|w\|^{2^*} \leq \varphi(t u_{\ge}) +  O(\ge^{N\frac{N-2}{2}}).
\]
The last estimate follows from the Young inequality
\[
\max_{s>0} \left( rs - \frac{s^p}{2} \right) = \frac{p-1}{p} r^{\frac{p}{p-1}}, \qquad p>1.
\]
We remark that $N(N-2)/(N+2) > 2$ since $N\geq 5$.
It now follows from Lemma \ref{lemma:calculus} and (\ref{eq:15}) that, for $\alpha$ and $\ge$ sufficiently
small, 
\begin{eqnarray*}
\max_{\substack{t > 0 \\ w \in Z}} \varphi(tu_{\ge} + w) &\leq& \frac{1}{N} \left( \frac{\int_{\Omega} |\nabla u_\ge|^2 - \lambda |u_{\ge}|^2 \, dx}{\left( \int_{\Omega} |u_{\ge}|^{2^*} |x|^\alpha \, dx  \right)^{2/2^*}}
\right)^{N/2} + O(\ge^{N\frac{N-2}{N+2}})\\
&<& \frac{1}{N} S^{\frac{N}{2}}.
\end{eqnarray*}
\noindent{\it Proof of Theorem \ref{th:main}}
Propositions \ref{prop:3.2} and \ref{prop:3.5} imply the existence of \(u \in \mathcal{N}\) such that
\(\varphi (u)=c\). In particular, \(D\varphi(u)=0\) on the tangent space \(T_u\mathcal{N}\). Since we have
shown that \(\mathcal{N}\) is a natural constraint, \(u\) is a free critical point of \(\varphi\).
\begin{remark}
When $\lambda>\lambda_1$, it is very easy to show that our solutions must change sign. Actually, just test (\ref{eq:1}) against $e_1$, and conclude that $u$ cannot have the same sign everywhere.
\end{remark}
In dimension \(N=4\), we can prove the following variant of Theorem \ref{th:main}.
\begin{theorem}
Assume \(N=4\) and that \(\lambda\) is not a Dirichlet eigenvalue of the Laplace operator. Then, for every $\alpha>0$ sufficiently small,
there exists (at least) a ground-state solution to problem (\ref{eq:1}).
\end{theorem}
The proof is achieved by an easy modification of the previous arguments. It suffices to take into accounts the
different asymptotic behavior of the instanton in dimension four.

\section{\bf Additional properties of ground-state solutions}
\def\theequation{2.\arabic{equation}}\makeatother
\setcounter{equation}{0}

As in \cite{MR2553063}, we can prove that ground-state solutions of (\ref{eq:1}) have more properties than
being just solutions. 

\begin{proposition}
Assume again that \(\lambda_m \leq \lambda < \lambda_{m+1}\). Then any point \(u \in \mathcal{N}\) such that
\(\varphi(u)=c\) is a critical point of \(\varphi\) with Morse index \(m+1\).
\end{proposition}
\noindent\proof
\(\mathcal{N}\) is a smooth manifold of codimension \(m+1\). With the notation introduced with Proposition \ref{prop:2.1},
we can write \(T_u \mathcal{N} = (DF(u))^{-1}(0)\). Since \(u\) minimizes \(\varphi\) on \(\mathcal{N}\), the hessian
of \(\varphi\) at \(u\) is positive definite on \(T_u\mathcal{N}\). We conclude that the Morse index of \(u\) is at most
\(m+1\). But the proof of Proposition \ref{prop:2.1} shows that this Morse index is at least \(m+1\), and the proof is complete.
Since \(\Omega\) is a radial domain, we might wonder if its symmetry is inherited by ground-state solutions.
We do not have a complete answer, as in the situation \(\alpha=0\) treated by \cite{MR2553063}.
However, we can still prove that ground-state solutions are \emph{foliated Schwarz symmetric} functions. We recall
the precise definition for the reader's sake.
\begin{definition}
A function \(u\) on a radial domain is \emph{foliated Schwarz symmetric} if there exists a unit vector \(p\in\mathbb{R}^N\)
such that \(u\) is a function of the distance from the origin and of the angle with the straight line along \(p\) only. In other words,
\(u(x)\) depends only on \(|x|\) and on \(\arccos \frac{x \cdot p}{|x|}\).
\end{definition}
When $\lambda_1 \leq \lambda < \lambda_2$, the constraint $\mathcal{N}$ is actually radially symmetric, by standard results about the symmetry of the first eigenfunction $e_1$. As the next results shows, in this situation we can gain more symmetry also for ground-state solutions.
\begin{proposition}
Let \(m=1\), i.e. \(\lambda_1 \leq \lambda < \lambda_2\). If \(u \in \mathcal{N}\) satisfies
\(\varphi(u)=\inf_{\mathcal{N}} \varphi = c\), then \(u\) is foliated Schwarz symmetric.
\end{proposition}
\noindent\proof
Under our assumptions, we remark that
\begin{align*}
\mathcal{N} &= \left\{ u \in H \setminus \{0\} \left| 
\begin{array}{ll}
&\int_\Omega |\nabla u|^2 - \lambda u^2 - |x|^\alpha u^{2^*} = 0 \\
&\int_\Omega \nabla u \cdot \nabla e_1 - \lambda u e_1 - |u|^{2^*-2}ue_1 =0
\end{array}
\right.
\right\} \\
&= \left\{ u \in H \setminus \{0\} \left| 
\begin{array}{ll}
&\int_\Omega |\nabla u|^2 - \lambda u^2 - |x|^\alpha u^{2^*} = 0 \\
&\int_\Omega (\lambda_1-\lambda)ue_1 - |u|^{2^*-2}ue_1 =0
\end{array}
\right.
\right\}
\end{align*}
We recall that \(e_1\), the first eigenfunction of the Laplace operator on \(\Omega\), is radially
symmetric and positive. Let \(u\) be as in the statement, and pick \(x_0 \in \Omega \setminus \{0\}\)
with
\[
u(x_0)=\max\{u(x)\mid x \in \overline{\Omega}, \ |x|=|x_0| \}.
\]
For \(p=x_0/|x_0|\), we define \(\mathcal{H}_p\) as the set of all closed halfspaces \(K\)
such that \(0 \in \partial K\) and \(p\) lies in the interior of \(K\).  For each \(K\in\mathcal{H}_p\),
there is a reflection map \(\sigma_K\) across \(K\). We need to prove (see \cite[Lemma 4.2]{MR2011581}) that
\begin{equation}\label{eq:5.1}
u(x) \geq u(\sigma_K(x)), \quad\text{for all \(x \in K \cap \Omega\)}.
\end{equation}
So, fix any $K\in\mathcal{H}_p$ and consider the polarization of $u$ with respect to $K$, defined by the formula
\[
u_K(x) = 
\begin{cases}
\max\{u(x),u(\sigma_K(x))\}, &\text{if $x \in \Omega \cap K$} \\
\min\{u(x),u(\sigma_K(x))\}, &\text{if $x \in\Omega \setminus K$}.
\end{cases}
\]
It is known that
\[
\int_\Omega |\nabla u_K|^2 = \int_\Omega |\nabla u|^2, \quad \int_\Omega |u_K|^q = \int_\Omega |u|^q
\]
for every $1 \leq q \leq +\infty$. Since $e_1$ is radially symmetric, we also have
\[
\int_\Omega u_K e_1 = \int_\Omega u e_1, \quad \int_\Omega |x|^\alpha |u_H|^{q-2}u_K e_1 = \int_\Omega |x|^\alpha |u|^{q-2}ue_1,
\]
for every $2 \leq q \leq +\infty$. We refer the interested reader to \cite[Section 2]{MR2177179}. As a consequence of these invariance properties, $u_K \in \mathcal{N}$ and $\varphi(u_K)=\varphi(u)=c$. Standard methods of elliptic regularity theory implies that both $u$ and $u_K$ are classical solutions of (\ref{eq:1}). Set $w=u_K-u$, and notice that $w \geq 0$ in $\Omega \cap K$; moreover, $w$ solves the Dirichlet problem
\begin{equation*}
\begin{cases}
-\Delta w = q(x)w &\text{in $(K\cap \Omega)^\circ$}\\
w = 0 &\text{on $\partial (K \cap \Omega)$},
\end{cases}
\end{equation*}
where $(K\cap \Omega)^\circ$ stands for the interior of $K\cap \Omega$. Here,
\[
q(x)=\lambda + (2^*-1) \int_0^1 |x|^\alpha |(1-s)u(x)+s u_K(x)|^{2^*-2}\, ds
\]
for every $x \in \Omega \cap K$. But $q \in L^\infty(\Omega \cap K)$, and the strong maximum principle tells us that either $w>0$ in $(\Omega \cap K)^\circ$ or $w$ identically vanishes. But $x_0 \in (\Omega \cap K)^\circ$ and $w(x_0)=u_K(x_0)-u(x_0) =0$, and thus $w =0$ everywhere. Hence $u_K=u$, and (\ref{eq:5.1}) is proved.
\begin{remark}
We observe that the previous proof is independent of the size of $\alpha$. Unlike \cite{MR2553063}, we are not able to exclude that $u$ is radially symmetric. Our equation contains the increasing weight $|\cdot|^\alpha$, and, as far as we know, there is no precise estimate for the Morse index of radially symmetric solutions of (\ref{eq:1}). See also \cite[Section 6.2]{MR2722502} for a recent survey on symmetry of solutions for similar equations.
\end{remark}

\section{\bf Final comments}
\def\theequation{1.\arabic{equation}}\makeatother
\setcounter{equation}{0}

Roughly speaking, the Dirichlet problem
\begin{equation} \label{limiting}
\begin{cases}
-\Delta u = \lambda u + |u|^{2^*-2}u &\text{in $\Omega$}\\
u = 0 &\text{on $\partial\Omega$}
\end{cases}
\end{equation}
is the \emph{limiting problem} for (\ref{eq:1}) as \(\alpha \to 0\). We have proved that many properties of 
this limiting problem pass on to (\ref{eq:1}) for small values of \(\alpha\). Although ours are pertubative results, it seems rather 
complicated to apply those methods developed in \cite{MR2186962}, since non-degeneracy of solutions to~(\ref{limiting}) is unknown.

On the other hand, when \(\lambda=0\), many authors studied the asymptotic properties of (\ref{eq:1}) as \(\alpha \to +\infty\): we refer to \cite{SSW} for seminal results. In our framework, we face a serious obstacle in (\ref{eq:15}). Indeed, one might try to push the spike \(x_\ell\) of the instanton towards \(\partial\Omega\), with a speed possibly related to \(\alpha\) as well. However, the denominator \((1-2\ell)^{2\alpha/2^*}\) behaves as an exponential function, whilst the numerator is a polynomial perturbation of the best Sobolev constant. We are therefore unable to treat this situation.

Let us try to explain this obstruction. By analogy with Theorem 3.4 of \cite{SW2003}, we may believe that the \emph{actual} limiting problem as \(\alpha \to +\infty\) is
\begin{equation*}
\begin{cases}
-\Delta V = \mathrm{e}^{x_1} |V|^{2^*-2}V &\text{in \(\mathbb{R}^N_{-}\)} \\
V=0 &\text{on \(\{x_1=0\}\)},
\end{cases}
\end{equation*}
where \(\mathbb{R}^N_{-} = \{x \in \mathbb{R}^N \mid x_1 <1\}\). Hence, the optimal level for compactness might be \emph{larger} than \(\frac{1}{N}S^{N/2}\). There would be room for existence of ground state solutions above \(\frac{1}{N}S^{N/2}\), but the instanton cannot suffice. Anyway, we do not have rigorous proofs of these ideas, yet.

\nocite{*}
\printbibliography

\enddocument